\documentclass[11pt]{article}
\usepackage{amsmath,amssymb,amscd,array}
\usepackage{amsmath}

\let\ssection=\section
\renewcommand{\section}{\setcounter{equation}{0}\ssection}

\def\d{\delta}

\def\om{\omega}
\def\r{\rho}
\def\a{\alpha}
\def\b{\beta}
\def\s{\sigma}
\def\vfi{\varphi}

\def\l{\lambda}
\def\m{\mu}

\def\implies{\Rightarrow}

\begin{document}

\frenchspacing

\def\d{\delta}
\def\g{\gamma}
\def\om{\omega}
\def\r{\rho}
\def\a{\alpha}
\def\b{\beta}
\def\s{\sigma}
\def\vfi{\varphi}
\def\l{\lambda}
\def\m{\mu}
\def\implies{\Rightarrow}

\oddsidemargin .1truein
\newtheorem{thm}{Theorem}[section]
\newtheorem{lem}[thm]{Lemma}
\newtheorem{cor}[thm]{Corollary}
\newtheorem{pro}[thm]{Proposition}
\newtheorem{ex}[thm]{Example}
\newtheorem{rmk}[thm]{Remark}
\newtheorem{rmks}[thm]{Remarks}
\newtheorem{defi}[thm]{Definition}
\title{Normalization of a nonlinear representation of a Lie algebra, regular on an abelian ideal}
\author{ Mabrouk BEN AMMAR}
\maketitle
\begin{abstract} We consider a nonlinear representation of a Lie algebra which is regular on an abelian ideal, we define a normal form which generalizes that defined in \cite{abs}.
\end{abstract}
\section{Introduction}
The study of a vector fields in a neighborhood of a point on a complex manifold is, of course, reduced to that of a vector fields $T$ in a neighborhood of the origin of $E=\mathbb{C}^n$.

If $T$ is regular at the origin of $E$, we know that there exists a coordinate system $(x_1,\dots,x_n)$ of $E$ in which the vector fields can be expressed as
$$
T=\frac{\partial}{\partial x_1},
$$
that is, there exists a local diffeomorphism $\phi$ of $E$ such that
$$
\phi(0)=0\quad\text{and}\quad \phi^\star(T)=\frac{\partial}{\partial x_1}.
$$
If $T$ is analytic then $\phi$ can be chosen to be analytic \cite{va}.

If $T$ is singular at the origin of $E$, then the situation is not so simple; generically, $T$ is linearizable, that is, there exists a local diffeomorphism $\phi$ of $E$ such that
$$
\phi(0)=0\quad\text{and}\quad \phi^\star(T)=\sum_{j,k}a_{jk}x_j\frac{\partial}{\partial x_k},\quad a_{jk}\in\mathbb{C}.
$$
But, the linearization is not always possible, so, we introduce the notion of normal form of a vector fields. Such a normal form, by construction, must enjoy the following properties:
\begin{itemize}
  \item [(a)] If $T$ is any vector fields then there exists a local diffeomorphism $\phi$ of $E$ such that $\phi(0)=0$ and $\phi^\star(T)$ is in normal form.
  \item [(b)] This normal form is unique, that is, if $\phi^\star_1(T)$ and $\phi^\star_2(T)$ are in normal form then $\phi^\star_1(T)=\phi^\star_2(T)$.
\end{itemize}
This problem has been studied extensively since Poincar\'e in formal, analytical or $C^infty$ contexts. The normal form of $T$ is then:
$$
T'=S^1+N
$$
wher $S^1$ is semisimple and $N$ is nilpotent satisfying
$$
[S^1,N]=0.
$$

Let us now consider the case, not of a single vector fields, but of a finite dimensional Lie algebra $\mathfrak{g}$ of vector fields, or, if we prefer, not the case of a germ of local actions on $E$, but rather the case of germs of local actions of a Lie group $G$  on $E$.

Let us consider the formal or analytic setting. If all vectors fields are singular at 0, then we are in the presence of a nonlinear representation $T$ of a Lie algebra $\mathfrak{g}$ in the sense of Flato,  Pinczon and Simon \cite{fps}:
$$
T: \mathfrak{g}\longrightarrow \mathcal{X}_0(E),\,X\longmapsto T_X\quad\text{such that}\quad [T_X,T_Y]=T_{[X,Y]},
$$
where $\mathcal{X}_0(E)$ is the space of vector fields which are singular at 0. In this case, we use the structure of the Lie algebra $\mathfrak{g}$ to precise the convenient normal forms. The notion of normal form given in \cite{abs} generalizes those known for the nilpotent and semisimple cases (\cite{abp}, \cite{b},\cite{gs},\cite{h}).
 Let $\mathfrak{r}$ be the solvable radical of $\mathfrak{g}$, we know that, by the Levi-Malcev decomposition theorem, $\mathfrak{g}$ can be decomposed
$$
(\mathcal{D})\qquad\mathfrak{g}=\mathfrak{r}\oplus \mathfrak{s}
$$
where $\mathfrak{s}$ is semi-simple. The nonlinear representation $T$ of the complex Lie algebra $\mathfrak{g}$ is said to be normal with respect the Levi-Malcev decomposition $(\mathcal{D})$ of $\mathfrak{g}$, if  its restriction $T|_{\mathfrak{s}}$ is linear and its restriction $T|_{\mathfrak{r}}$ to $\mathfrak{r}$ is in the following form
$$
T|_{\mathfrak{r}}=D^1+N^1+\sum_{k\geq2}T^k|_{\mathfrak{r}},
$$
where the linear part $D^1+N^1$ of $T|_{\mathfrak{r}}$ is such that $D^1$ is diagonal with coefficients $\mu_1,\,\dots,\,\mu_n$, (the $\mu_i$ are elements of $\mathfrak{r}^*$) and $N^1$ is strictly upper triangular. Moreover, $T^k$ has the following form in the coordinates $x_i$ of $E$
$$
T_X^k=\sum_{i=1}^n\sum_{|\alpha|=k}\Lambda^\alpha_i(X)x_1^{\alpha_1}\cdots x_n^{\alpha_n} \frac{\partial}{\partial x_i},\qquad X\in\mathfrak{r},
$$
where any coefficient $\Lambda^\alpha_i$ can be nonzero only if it is resonant, that means, the corresponding linear form $\mu_\a^i\in\mathfrak{r}^*$ defined by
$$
\mu_\a^i=\sum_{r=1}^n \alpha_r\mu_r-\mu_i,
$$
is a root of $\mathfrak{r}$, (eigenvalue of the adjoint representation of $\mathfrak{r}$).
The representation $T$ is said to be normalizable with respect $(\mathcal{D})$  if it is equivalent to a normal one with respect  $(\mathcal{D})$.

In this paper, we study the situation where the $T_X$ are not all in $\mathcal{X}_0(E)$. More precisely, we introduce the notion of nonlinear representation $(T,E)$ of $\mathfrak{g}$ in $E$ such that
$$
T_X=T^0_X+\sum_{k\geq1}T^k_X\in\mathcal{X}(E),\qquad X\in\mathfrak{g},
$$
where $\mathcal{X}(E)$ is the space of formal vector fields on $E$ (not necessarily vanishing at 0), in particular, $T^0_X\in E$. This situation seems to be more delicate in its generality, but here we treat the case where the Lie algebra $\mathfrak{g}$ can be written $$\mathfrak{g} = \mathfrak{g}_0 + \mathfrak{m},$$ where $\mathfrak{m}$ is an abelian ideal of $\mathfrak{g}$ and $\mathfrak{g}_0$ is a subalgebra of $\mathfrak{g}$. The Lie algebras of groups of symmetries of simple physical systems are often of this type (Galil\'ee group, Poincar\'e group, ...). Moreover, we assume that
 $$
 T^0_X=0 \quad\text{if}\quad X\in\mathfrak{g}_0 \quad\text{and}\quad T^0_X\neq 0 \quad\text{if}\quad X\in\mathfrak{m}\setminus\{0\}.
 $$
We also assume that $T$ is analytic and we will prove in this case that we can always put the representation $T$ in the following normal form: there exists a coordinate system $(x_1,\dots,x_p,y_1,\dots, y_q)$ of $E$ such that:
\begin{itemize}
  \item [(a)] For all $X\in\mathfrak{m}$, $T_X$ can be expressed as
$$
T_X=\sum_{i=1}^pa_i\frac{\partial}{\partial x_i},
$$
 where $a_i\in\mathfrak{m}^*$.

  \item [(b)] $T|_{\mathfrak{g}_0}$ is normal au sense de \cite{abs}, moreover, for all $X\in \mathfrak{g} _0$ and $k\geq2$, $T^k_X$ has the following form
$$
T^k_X=\sum_{i=1}^p\sum_{|\alpha|=k}\Lambda^\alpha_i(X)y_1^{\alpha_1}\cdots y_q^{\alpha_q} \frac{\partial}{\partial x_i}+\sum_{i=1}^q\sum_{|\alpha|=k}\Gamma^\alpha_i(X)y_1^{\alpha_1}\cdots y_q^{\alpha_q} \frac{\partial}{\partial y_i}
$$
\end{itemize}

 \section{Notations and definitions}
 Let $E$ be a complex vector space with dimension $n$. The space of symmetric $k$-linear applications from $E\times\cdots\times E$ to $E$ is identified with the space
$\mathrm{L}(\otimes_s^kE,E)$ of linear maps from $\otimes_s^kE$ to $E$, where $\otimes_s^kE$ is the space of symmetric $k$-tensors on $E$. Denote by $\mathrm{L}(E)$ the space $\mathrm{L}(E,E)$.
Let $(e_1,\dots,e_n)$ be a basis of $E$. For $\alpha\in\mathbb{N}^n$ such that $|\alpha|=k$ and $i\in\{1,\dots,n\}$
 we define $e_\alpha^i\in \mathrm{L}(\otimes_s^kE,E)$ by:
$$
e^\a_i\circ\sigma_k(\otimes^{\beta_1}e_1\cdots\otimes^{\beta_n}e_n )=\delta_{\alpha_1}^{\beta_1}\cdots\delta_{\alpha_n}^{\beta_n}e_i,
$$
where $\delta$ is the Kronecker symbol, and $\sigma_k$ is the symmetrization operator from $\otimes^kE$ to $\otimes^k_sE$ defined by
$$
\s_k(v_1,\dots,v_k)=\frac{1}{k!}\sum_{\s\in\mathcal{S}_k}(v_{\s(1)},\dots,v_{\s(k)}).
$$
Obviously, $(e^\a_i)$ is a basis of $\mathrm{L}(\otimes_s^kE,E)$. Let $\mathcal{X}(E)$ (respectively $\mathcal{X}_0(E)$) be the set of formal power series (or formal vector fields)
  $$
  T=\sum_{k=1}^\infty T^k,\quad\text{(respectively~~~} T=\sum_{k=0}^\infty T^k)
  $$
  where $T^k\in\mathrm{L}(\otimes_s^kE,E)$ (with $\mathrm{L}(\otimes_s^0E,E)=E$). Therefore, we can write
  $$
  T^k=\sum_{|\alpha|=k}\sum_{i=1}^ne^\a_i.
  $$
  Of course, any analytical vector fields $X$ on $E$ which is singular in 0 admits a Taylor expansion:
  $$
 \sum_{i=1}^n \sum_{k=1}^\infty\sum_{|\a|=k}\Lambda^\a_ix_1^{\alpha_1}\cdots x_n^{\a_n}\frac{\partial}{\partial x_i}=\sum_{\a, i}\Lambda^\a_ix^\a\partial_i.
  $$
  It can be identified with the formal vector fields:
  $$
  \sum_{i=1}^n \sum_{k=1}^\infty\sum_{|\a|=k}\Lambda^\a_ie^\a_i.
  $$
  If $|\a|=0$ we agree that $e^\a_i=e_i$.

Of course, $\mathcal{X}(E)$ can be endowed with a Lie algebra structure defined, for $e^\a_i\in \mathrm{L}(\otimes^{|\a|}_s E,E)$ and $e^\b_j\in \mathrm{L}(\otimes^{|\b|}_s E,E)$, by
$$
[e^\a_i,e^\b_j]= e^\a_i\star e^\b_j- e^\b_j\star e^\a_i\in\mathrm{L}(\otimes^{|\a|+|\b|-1}_s E,E),
$$
where
$$
e^\a_i\star e^\b_j=\\b_ie^{\a+\b-1_i}_j,\qquad 1_i=(0,\dots,0,1,0,\dots,0), (1 \text{ in } i\text{th place}).
$$
\begin{defi}
A formal vector fields $T$ in $\mathcal{X}(E)$ is said to be analytic al if the power series
$$
\sum_{k\geq1}T^k(\otimes^kv)=T_X(v),\quad v\in E,
$$
converges in a neighborhood of the origin of $E$.
\end{defi}
\begin{defi}
Two formal vector fields $T$ and $T'$ are said to be equivalent if there exists an element $\phi=\sum_{k\geq1}\phi^k$ of $\mathcal{X}_0(E)$ such that $\phi^1$ is invertible and
$$
\psi\star T=T'\circ \psi
$$
where
$$
T'\circ\phi=\sum_{k\geq1}\sum_{j=1}^kT'^j\circ\left(\sum_{i_1+\cdots+i_j=k}\phi^{i_1}\otimes\cdots\otimes\psi^{i_j}\right)\circ\sigma_k.
$$
$T$ and $T'$ are said to be analytically equivalent if $\phi$ is analytic.
\end{defi}
\begin{rmk}
For the composition law $\circ$, the map $\phi$ is invertible if and only if $\phi^1$ is an automorphism of $E$.
\end{rmk}
\begin{defi}
Let $\mathfrak{g}$ be a Lie algebra with finite dimension over $\mathbb{C}$. A nonlinear (formal) representation $(T,E)$ of $\mathfrak{g}$ in $E$ is a linear map:
$$
T:\mathfrak{g}\rightarrow\mathcal{X}(E),\quad X\mapsto T_X,
$$
such that
$$
[T_X,T_Y]=T_{[X,Y]},\quad X,\,Y\in\mathfrak{g}.
$$
\end{defi}
\begin{defi}
We say that two representations $(T,E)$ and $(T',E)$ are equivalent if there exists $\phi\in\mathcal{X}_0(E)$ such that $\phi_1$ is invertible and
$$
\phi\star T_X=T'_X\circ\phi,\quad X\in\mathfrak{g}.
$$
\end{defi}
\begin{defi}
A nonlinear formal representation $T$ of $\mathfrak{g}$ is said to be analytic if the power series
$$
\sum_{k\geq1}T^k_X(\otimes^ke)=T_X(v),\qquad v\in E
$$
converges in a neighborhood of 0 in $E$, for all $X\in\mathfrak{g}$.
\end{defi}

In the following, we consider a nonlinear analytic representation $(T,E)$ of a complex finite dimensional Lie algebra $\mathfrak{g}$ in $E=\mathbb{C}^n$, such that
$$
\mathfrak{g} = \mathfrak{g}_0 + \mathfrak{m},
$$
 where $\mathfrak{m}$ is an abelian ideal of $\mathfrak{g}$ and $\mathfrak{g}_0$ is a subalgebra of $\mathfrak{g}$. Moreover, we assume that
 $$
 T^0_X=0 \quad\text{if}\quad X\in\mathfrak{g}_0 \quad\text{and}\quad T^0_X\neq 0 \quad\text{if}\quad X\in\mathfrak{m}\setminus\{0\},
 $$
 and then we normalize $T$ step by step.
 \section{Normalization of $T|_\mathfrak{m}$}
 \begin{pro}
 There exists a coordinate system $(x_1,\dots,x_p,y_1,\dots, y_q)$, $p+q=n$, of $E$ such that, for all $X\in\mathfrak{m}$, $T_X$ can be expressed as
$$
T_X=\sum_{i=1}^pa_i\frac{\partial}{\partial x_i},
$$
where $a_i\in\mathfrak{m}^*$.
 \end{pro}
 {\it Proof.} Since the representation $(T,E)$ is analytic, then, there exists a neighborhood $\mathcal{U}$ of the origin of $E$ such that
 $$
 T_x(v)\in E,
 $$
  for any $v$ in $\mathcal{U}$ and $X$ in $\mathfrak{g}$. So, for any $v$ in $\mathcal{U}$, we consider the subspace $E_v$ spanned by the vectors $T_X(v)$, where $X$ browses $\mathfrak{m}$. If $(X_1,\dots,X_p)$ is a basis of $\mathfrak{m}$, then  the $T_{X_i}(v)$ are generators of $E_v$, thus, $\mathrm{dim}E_v\leq p$. But, the system $(T_{X_1}(0),\dots,T_{X_p}(0))$ is independent, indeed, the condition
  $$
  \a_1T_{X_1}(0)+\cdots+\a_pT_{X_p}(0)=T^0_{\a_1{X_1}(0)+\cdots+\a_p{X_p}}=0
  $$
  implies that $\a_1{X_1}(0)+\cdots+\a_p{X_p}=0$ since $T^0$ does not vanish on $\mathfrak{m}\setminus\{0\}$. Thus,  $\mathrm{dim}E_0= p$. The determinant map is a continuous map, then there exists a neighborhood $\mathcal{V}\subset\mathcal{U}$ such that $\mathrm{dim}E_v\geq p$, for all $v\in\mathcal{V}$. Therefore, $\mathrm{dim}E_v=p$, for all $v\in\mathcal{V}$.

  Thus, the map
  $$
  v\mapsto E_v
  $$
  is an involutive integrable distribution, therefore, the Frobinius theorem ensures the existence of an analytic coordinate system $(x_1,\dots,x_p,y_1,\dots,y_q)$ of $E$ such that the $\frac{\partial}{\partial x_i}$, $i=1,\dots,p$, generate this distribution.

  \hfill$\Box$

 \section{Normalization of $T|_{\mathfrak{g}_0}$}
 Now, we begin the second step to normalize the representation $(T,E)$. We know that, for any $X$ in $\mathfrak{g}_0$, we have
  $$
  T_X=\sum_{k\geq1}T^k_X.
  $$
  We consider the coordinate system $(x_1,\dots,x_p,y_1,\dots,y_q)$ of $E$ defined in the previous section and we prove the following results
 \begin{pro}
 For any $X$ in $\mathfrak{g}_0$,  $T^1_X$ has the following form
 $$
 T^1_X=\sum_{i,j}a_{ij}(X)y_j\frac{\partial}{\partial x_i}+\sum_{i,j}b_{ij}(X)y_j\frac{\partial}{\partial y_i}+\sum_{i,j}c_{ij}(X)x_j\frac{\partial}{\partial x_i}
 $$
 and for $k\geq2$, $T^k_X$ has the following form
 $$
 T^k_X=\sum_{i,|\a|=k}A_{i}^\a(X)y^\a\frac{\partial}{\partial x_i}+\sum_{i|\a|=k}B_{i}^\a(X)y^\a\frac{\partial}{\partial y_i}.
 $$
 \end{pro}
 {\it Proof.} Consider $(X_1,\dots,X_p)\in\mathfrak{m}^p$ such that
 $$
 T_{X_i}=\frac{\partial}{\partial x_i}.
 $$
 Since $\mathfrak{m}$ is an ideal of $\mathfrak{m}$, then $[\mathfrak{m},\mathfrak{g}_0]\subset\mathfrak{m}$, and therefore, for any $X$ in $\mathfrak{g}_0$ and $j=1,\dots,p$, we have
 $$
 T_{[X_j,X]}=\left[\frac{\partial}{\partial x_j},T_X\right]=\left[\frac{\partial}{\partial x_j},T^1_X\right]=\sum_{i=1}^p\a_i\frac{\partial}{\partial x_i},
 $$
 for some $\a_i\in\mathbb{C}$. In particular, for $k\geq2$, we have
 $$
 \left[\frac{\partial}{\partial x_j},T^k_X\right]=0,
 $$
 and then the result follows.

 \hfill$\Box$

Now, for any $X$ in $\mathfrak{g}_0$, we define
  $$
  \begin{array}{llllll}
    A_X&=&\sum_{i,\a}A_{i}^\a(X)y^\a\frac{\partial}{\partial x_i}=\sum_{k\geq1}A^k_X,\quad (A_i^{1_j}=a_{ij}),\\[8pt]
    B_X&=&\sum_{i,\a}B_{i}^\a(X)y^\a\frac{\partial}{\partial x_i}+\sum_{i,j}c_{ij}(X)x_j\frac{\partial}{\partial x_i}=\sum_{k\geq1}B^k_X,\\[8pt]
    H_X&=&\sum_{i,j}b_{ij}(X)y_j\frac{\partial}{\partial y_i},\\[8pt]
    K_X&=&\sum_{i,j}c_{ij}(X)x_j\frac{\partial}{\partial x_i},
  \end{array}
  $$
    and we consider the subspaces $E_1$ and $E_2$ corresponding, respectively, to coordinate systems $(x_1,\dots,x_p)$ and $(y_1,\dots,y_q)$.
  \begin{pro}
  i) $(H,E_1)$ and $(K,E_2)$ are two linear representations of $\mathfrak{g}_0$.

  ii) For any $X$ in $\mathfrak{g}_0$ we have
  $$
  \begin{array}{lll}
  A_{[X,Y]}=[A_X,B_Y]+[B_X,A_Y],\\[8pt]
  B_{[X,Y]}=[B_X,B_Y].
  \end{array}
  $$
  Thus, $(B,E)$ is a nonlinear representation of $\mathfrak{g}_0$ in $E$.
  \end{pro}

  {\it Proof.} i) Obvious.

  ii) For any $X$ in $\mathfrak{g}_0$, we have
  $$
  T_X=A_X+B_X,
  $$
then we have
$$
T_{[X,Y]}=A_{[X,Y]}+B_{[X,Y]}=[A_X+B_X,A_Y+B_Y]=[A_X,B_Y]+[B_X,A_Y]+[B_X,B_Y].
$$
We easily check that
$$
  \begin{array}{lll}
  A_{[X,Y]}=[A_X,B_Y]+[B_X,A_Y],\\[8pt]
  B_{[X,Y]}=[B_X,B_Y].
  \end{array}
  $$

  \hfill$\Box$

Now, we consider a Levi-Malcev decomposition of $\mathfrak{g}_0$:
$$
(\mathcal{D})\qquad\mathfrak{g}=\mathfrak{r}\oplus \mathfrak{s}
$$
where $\mathfrak{s}$ is semi-simple and $\mathfrak{r}$ is the solvable radical of $\mathfrak{g}$. We triangularize simultaneously the $H_X$ and $K_X$, where $X$ browses $\mathfrak{r}$ (Lie theorem \cite{j}). Thus, we define the linear forms: $\mu_1,\dots,\mu_p\in\mathfrak{r}^*$ and $\nu_1,\dots,\nu_q\in\mathfrak{r}^*$. The $\mu_1(X),\dots,\mu_p(X),\nu_1(X),\dots,\nu_q(X)$ are the eigenvalues of $B^1(X)$, indeed, $B^1_X=H_X+K_X$. They are also  the eigenvalues of $T^1(X)$, since $A^1_X\in\mathrm{L}(E_2,E_1)$. Thus, the simultaneous triangulation of $H$ and $K$ leads to that of $T^1$ by retaining the first $p$ components in $E_1$ and the $q$ other in $E_2$.

Let us denote by $\l_1,\dots,\l_n$ the roots of $T^1|_\mathfrak{r}$, and consider the elements $\l^\a_j$ of $\mathfrak{r}^*$ defined by
$$
\l^\a_j(X)=\sum_{i=1}^n\a_i\l_i(X)-\l_j(X),
$$
where $j=1,\dots,n$ and $\a\in\mathbb{N}$.
\begin{defi}
We say that $(\a,j)$ is resonant if $\l^\a_j$ is a root of $\mathfrak{r}$.
\end{defi}
\begin{defi}
An element  $X_0$ in $\mathfrak{r}$ is said to be vector resonance if, for any $\a$, for any $i$ and for any $j$,
$$
\l^\a_j(X_0)=v_i(X_0)\Rightarrow \l^\a_j=v_i,
$$
where the $v_i$ are the roots of $\mathfrak{r}$.
\end{defi}
For the two following lemmas see, for instance, \cite{abp} and \cite{abs}.
\begin{lem}
The set of vector resonance is dense in $\mathfrak{r}$.
\end{lem}
\begin{lem}
There exists a resonance vector in $\mathfrak{r}$ such that $[X_0,\mathfrak{s}]=0$.
\end{lem}
The resonant pairs $(\a, j)$ appearing here are of two types:
\begin{itemize}
  \item [(1)] $\sum_{i=1}^p\a_i\nu_i-\mu_j$ is a root of $\mathfrak{r}$.
  \item [(2)] $\sum_{i=1}^p\a_i\nu_i-\nu_j$ is a root of $\mathfrak{r}$.
\end{itemize}
Let us denote by
$$
R=\{(\a,j)~|~(\a,j)\text{ resonant of type } (1)\}
$$
and by
$$
R'=\{(\a,j)~|~(\a,j)\text{ resonant of type } (2)\}
$$
\begin{thm}
$T|_\mathfrak{r}$ is normalizable in the sense of \cite{abs}, that means, there exists an analytic operator $\phi=\sum \phi^k$ in $\mathcal{X}_0(E)$, such that $\phi^1$ is invertible and, for any $X$ in $\mathfrak{r}$, $\phi\star T_X\circ\phi^{-1}$ is in the form:
$$
\sum_{(\a,i)\in R}\Lambda^\a_i(X)y^\a\frac{\partial}{\partial x_i}+\sum_{(\a,i)\in R'}\Gamma^\a_i(X)y^\a\frac{\partial}{\partial y_i}.
$$
\end{thm}
{\it Proof.}
Let $X_0$ be a resonance vector of $\mathfrak{r}$ such that $[X_0,\mathfrak{s}]=0$. It is well known that $T_{X_0}$ is analytically normalizable: there exists an analytic operator $\phi=\sum \phi^k$ in $\mathcal{X}_0(E)$, such that $\phi^1$ is invertible and, for any $X$ in $\mathfrak{r}$,
$$
\phi\star T_{X_0}\circ\phi^{-1}=T'_{X_0},
$$
where $T'_{X_0}=S^1_{X_0}+N_{X_0}$ with
$$
 S^1_{X_0}=\sum_{i=1}^p\mu_{i}(X_0)x_i\frac{\partial}{\partial x_i}+\sum_{i=1}^q\nu_{i}(X_0)y_i\frac{\partial}{\partial y_i},
 $$
 and
 $$
[S^1_{X_0},N_{X_0}]=0.
$$
Therefore, $N_{X_0}$ has the following form
$$
N_{X_0}=\sum_{(\a,i)\in R_0}\Lambda^\a_i(X_0)y^\a\frac{\partial}{\partial x_i}+\sum_{(\a,i)\in R'_0}\Gamma^\a_i(X_0)y^\a\frac{\partial}{\partial y_i},
$$
where
$$
R_0=\{(\a,i)~|~\sum_{j=1}^q\a_j\nu_j-\mu_i=0\}\quad\text{and}\quad R'_0=\{(\a,i)~|~\sum_{j=1}^q\a_j\nu_j-\nu_i=0\}.
$$
The operator $\phi$ normalizes $T|_\mathfrak{r}$ (see \cite{abs}).

\hfill$\Box$
\begin{pro}
The normalization operator $\phi$ of $T|_\mathfrak{r}$ leaves invariant $T|_\mathfrak{m}$.
\end{pro}
{\it Proof.} To prove that $\phi$ leaves invariant $T|_\mathfrak{m}$ we will prove that
$$
\phi=\cdots(I+W_k)\circ\cdots\circ(I+W_1),
$$
where $I$ is the identity of $E$ and $W_k\in\mathrm{L}(\otimes^k_s E,E)$.

(a) $B^1_{X_0}=H_{X_0}+K_{X_0}$ being decomposed into a semisimpe part $S^1_{X_0}$ and a nilpotent part, therefore, to reduce $T^1_{X_0}$ it suffices to reduce $A^1_{X_0}$. We decompose $A^1_{X_0}$:
$$
A^1_{X_0}=A^1_{0X_0}+A^1_{1X_0},
$$
where
$$
A^1_{0X_0}=\sum_{(i,j), \nu_j=\mu_i}a_{ij}(X_0)y_i\frac{\partial}{\partial x_i}\in\mathrm{ker}~\mathrm{ad} S^1_{X_0}
$$
and
$$
A^1_{1X_0}=\sum_{(i,j), \nu_j\neq\mu_i}a_{ij}(X_0)y_i\frac{\partial}{\partial x_i}\in\mathrm{Im}~\mathrm{ad} S^1_{X_0}.
$$
We reduce $A^1_{X_0}$ by removing $A^1_{1X_0}$. Indeed, for any $X$ and $Y$ in $\mathfrak{r}$,
$$
A^1_{[X,Y]}=[A^1_X,B^1_Y]+[B^1_X,A^1_Y],
$$
therefore, there exists $W^1$ in $\mathrm{L}(E)$ and a 1-cocycle $V^1$ such that
$$
A^1_X=[B^1_X,W_1]+V^1.
$$
In particular, $A^1_{X_0}=[B^1_{X_0},W_1]+V^1$, then we can choose, for instance,
$$
W_1=(\mathrm{ad}~B^1_{X_0}(A^1_{1X_0})\quad\text{and}\quad V^1_{X_0}=A^1_{0X_0}.
$$
The stability of the eigenvector subspaces of $\mathrm{ad}~S^1_{X_0}$ by $B^1_{X_0}$ proves that $W^1\in \mathrm{L}(E_2,E_1)$, therefore, $(I+W_1)$ is invertible and
$$
(I+W_1)\star T^1_{X_0}\circ(I+W^1)^{-1}=B^1_{X_0}+A^1_{0X_0},~~~~~~~~~~~~~~~~~~~~~~~~(1)
$$
$$
(I+W_1)\star \frac{\partial}{\partial x_i}\circ(I+W_1)^{-1}=\frac{\partial}{\partial x_i}, \quad i=1,\dots,p.~~~~~~~~~~~~~~(2)
$$
$T|_\mathfrak{m}$ is then stable by $(i+W_1)$.

b) Now, we assume that $T_{X_0}$ is normalized up to order $k-1$, we write
$$
T_{X_0}=S^1_{X_0}+\sum_{1\leq j<k}N^j_{X_0}+T^k_{X_0}+\sum_{j>k}T^j_{X_0}
$$
and we decompose $T^k_{X_0}$:
$$
T^k_{X_0}=T^k_{0X_0}+T^k_{1X_0},
$$
where
$$
T^k_{0X_0}\in\mathrm{ker}~\mathrm{ad} S^1_{X_0}\quad\text{and}\quad T^k_{1X_0}\in\mathrm{Im}~\mathrm{ad} S^1_{X_0}.
$$
Therefore, $T^k_{0X_0}$ has the following form
$$
T^k_{0X_0}=N^k_{X_0}=\sum_{(\a,i)\in R_0}\Lambda^\a_i(X_0)y^\a\frac{\partial}{\partial x_i}+\sum_{(\a,i)\in R'_0}\Gamma^\a_i(X_0)y^\a\frac{\partial}{\partial y_i}.
$$
We seek an operator $W_k$ in $\mathrm{L}(\otimes^k_sE,E)$ such that
$$
(1+W_k)\star T_{X_0}=(S^1_{X_0}+\sum_{1\leq j\leq k}N^j_{X_0}+\sum_{j>k}T^j_{X_0})\circ(1+W_k).
$$
We can choose
$$
W_k=(\mathrm{ad}~T^1_{X_0})^{-1}(T^1_{1X_0}).
$$
The eigenvector subspaces of $\mathrm{ad}~S^1_{X_0}$ are stable by $\mathrm{ad}~T^1_{X_0}$, therefore, $W_k$
has the following form:
$$
W_k=\sum_{(\a,i)\in R_0}W^\a_iy^\a\frac{\partial}{\partial x_i}+\sum_{(\a,i)\in R'_0}W'^\a_iy^\a\frac{\partial}{\partial y_i}.
$$
Thus,
$$
(I+W_k)\star\frac{\partial}{\partial x_i}\circ(I+W_k)^{-1}=\frac{\partial}{\partial x_i},\quad i=1,\dots,p.
$$
We proceed by induction, so, we construct
$$
\phi=\cdots(I+W_k)\circ\cdots\circ(I+W_1)
$$
which normalizes $T|_\mathfrak{r}$ and leaves invariant $T|_\mathfrak{m}$.

\section{Linearization of $T|_\mathfrak{s}$}
Let us now consider the representation $(T',E)$ of $\mathfrak{g}$ in $E$, defined on $\mathfrak{g}$ by
$$
T'_X=\phi\star T_X\circ\phi^{-1}.
$$
The representation being $T'|_\mathfrak{r}$, as in \cite{abs}, we construct an analytic invertible operator $\psi\in\mathcal{X}_0(E)$, which linearizes $T'|_\mathfrak{s}$ and leaves $T'|_\mathfrak{r}$ in normal form. This is possible through the choice of $X_0$ switching with $\mathfrak{s}$: $[X_0,\mathfrak{s}]=0$.

On the other hand, by construction, the components of $\psi$ are independent of the coordinates $x_1,\dots,x_p$, then $\psi$ leaves invariant $T'|_\mathfrak{m}=T|_\mathfrak{m}$.
\section{Recapitulation}
Consider a nonlinear analytic representation $(T,E)$ of a complex finite dimensional Lie algebra $\mathfrak{g}$ in $E=\mathbb{C}^n$, such that
$$
(\mathcal{D}_1)\qquad\mathfrak{g} = \mathfrak{g}_0 + \mathfrak{m},
$$
 where $\mathfrak{m}$ is a $p$-dimensional abelian ideal of $\mathfrak{g}$ and $\mathfrak{g}_0$ is a subalgebra of $\mathfrak{g}$. Moreover, we assume that
 $$
 T^0_X=0 \quad\text{if}\quad X\in\mathfrak{g}_0 \quad\text{and}\quad T^0_X\neq 0 \quad\text{if}\quad X\in\mathfrak{m}\setminus\{0\}.
 $$
We consider a Levi-Malcev decomposition of $\mathfrak{g}_0$:
$$
(\mathcal{D}_2)\qquad\mathfrak{g}=\mathfrak{r}\oplus \mathfrak{s}
$$
where $\mathfrak{s}$ is semi-simple and $\mathfrak{r}$ is the solvable radical of $\mathfrak{g}$. Then, we have
\begin{thm} The representation $(T,E)$ of $\mathfrak{g}$ can be normalized with respect the decompositions $(\mathcal{D}_1)$ and $(\mathcal{D}_2)$. That means, there exists a coordinate system $(x_1,\dots,x_p,y_1,\dots, y_q)$ of $E$ such that
\begin{itemize}
  \item [a)] For all $X\in\mathfrak{m}$, $T_X$ can be expressed as
$$
T_X=\sum_{i=1}^pa_i\frac{\partial}{\partial x_i},
$$
where $a_i\in\mathfrak{m}^*$.
  \item [b)] $T|_\mathfrak{s}$ is a linear representation of $\mathfrak{s}$ in $E$.
  \item [c)] For any $X$ in $\mathfrak{r}$, $T_X$ is in the form:
$$
T_X=\sum_{(\a,i)\in R}\Lambda^\a_i(X)y^\a\frac{\partial}{\partial x_i}+\sum_{(\a,i)\in R'}\Gamma^\a_i(X)y^\a\frac{\partial}{\partial y_i},
$$
where $R_0=\{(\a,i)~|~\sum_{j=1}^q\a_j\nu_j-\mu_i=0\}$ and $ R'_0=\{(\a,i)~|~\sum_{j=1}^q\a_j\nu_j-\nu_i=0\}$.
\end{itemize}
\end{thm}


\end{document}